\documentclass[a4paper,12pt,intlimits,oneside]{amsart}

\usepackage{enumerate}
\usepackage{latexsym,amssymb}

\newcommand{\comment}[1]{}
\numberwithin{equation}{section}

\newcommand{\epf}{ $\Box$\medskip}

\newcounter{rea}
\newcounter{rej}
\newcounter{res}
\newtheorem{thm}{Theorem}[section]
\newtheorem{prop}[thm]{Proposition}
\newtheorem{cor}[thm]{Corollary}
\newtheorem{lem}[thm]{Lemma}
\newtheorem{defn}[thm]{Definition}

\begin{document}

\title[]{A weighted inequality for potential type operators}

\author[A. Adama]{A\"issata Adama}
\address{DER Math\'ematiques et informatique Universit\'e de Bamako, BPE 3206 Bamako-Mali}
\email{{\tt aissatalawdi@yahoo.fr}}
\author[J. Feuto]{Justin Feuto}
\address{Laboratoire de Math\'ematiques Fondamentales, UFR Math\'ematiques et Informatique, Universit\'e F\'elix Houphou\"et-Boigny (Abidjan), 22 B.P 1194 Abidjan 22. C\^ote d'Ivoire}
\email{{\tt justfeuto@yahoo.fr}}
\author[I. Fofana]{Ibrahim Fofana}
\address{Laboratoire de Math\'ematiques Fondamentales, UFR Math\'ematiques et Informatique, Universit\'e F\'elix Houphou\"et-Boigny (Abidjan), 22 BP 582 Abidjan 22, C\^ote d'Ivoire} 
\email{{\tt  fofana\_ib\_math\_ab@yahoo.fr}}

\thanks{}
\begin{abstract}
We establish a weighted inequality for fractional maximal and convolution type operators, between weak Lebesgue spaces and  Wiener amalgam type spaces on $ \mathbb R $ endowed with a measure which needs not to be doubling.
\end{abstract}
\maketitle

\section{Introduction}
Let $d$ be a positive integer and $\mu$ a nonegative Radon measure on $\mathbb R^d$ satisfying
$$\mu(\partial Q)=0,\ Q\text{ cube}.$$
By a cube, we always mean a bounded cube in $\mathbb R^d$ with sides parallel to the coordinate axes.

Let $1\leq q<\beta\leq\infty$. The fractional maximal operator $\mathfrak m^{\mu}_{q,\beta}$ and the centered one $\mathfrak m^{\mu,c}_{q,\beta}$ are defined as follows: for any element $f$ of $L^{1}_{\mathrm{loc}}(\mu)$ and any point $x\in\mathbb R^d$
\begin{equation}
\mathfrak m^{\mu}_{q,\beta}f(x)=\sup_{Q\ni x}\mu(Q)^{\frac{1}{\beta}-\frac{1}{q}}\left(\int_{Q}\left|f(x)\right|^qd\mu(x)\right)^{\frac{1}{q}},\label{opmax}
\end{equation}
where the supremum is taken over all cubes $Q$ containing $x$.
\begin{equation}
\mathfrak m^{\mu,c}_{q,\beta}f(x)=\sup_{Q\in\mathcal Q(x)}\mu(Q)^{\frac{1}{\beta}-\frac{1}{q}}\left(\int_{Q}\left|f(x)\right|^qd\mu(x)\right)^{\frac{1}{q}}
\end{equation}
where $\mathcal Q(x)=\left\{Q \text{ cube centered at }x\right\}$.

These operators plays an important role in harmonic analysis and partial differential equations theory. Their boundedness between weighted and non-weighted Lebesgue spaces have been studied in depth (see \cite{MW},\cite{OP},\cite{PW},\cite{SS} and the references therein). The results are mostly expressed in terms of weights belonging to the so-called Muckenhoupt classes $\mathcal A_\alpha(\mu)$ defined as follows.
\begin{defn}
Let $1\leq \alpha<\infty$ and $w$ a weight (a nonnegative locally $\mu$-integrable function) on $\mathbb R^d$.
\begin{enumerate}
\item[a)] When $1<\alpha<\infty$, we says that $w\in\mathcal A_\alpha$ if there exists a constant $C>0$ such that for every cube $Q$,
$$\left(\frac{1}{\mu(Q)}\int_{Q}w(x)d\mu(x)\right)\left(\frac{1}{\mu(Q)}\int_{Q}w(x)^{1-p'}d\mu(x)\right)^{p-1}\leq C$$
where $p'$ is the conjugate exponent of $p$.
\item[b)] We say that $w\in\mathcal A_1(\mu)$  if there exists a constant $C>0$ such that for every cube $Q$,
$$\frac{1}{\mu(Q)}\int_{Q}w(x)d\mu(x)\leq C\mathrm{ess}\inf_{x\in Q}w(x).$$
\end{enumerate}
\end{defn}
A well known result of Muckenhoupt and Wheeden reads as.
\begin{thm}[Theorem 2 \cite{MW}]\label{MuWee}Let $1\leq \alpha<\beta$, $\frac{1}{s}=\frac{1}{\alpha}-\frac{1}{\beta}$. If $\mu$ is the Lebesgue measure on $\mathbb R^d$ and $w$ a weight belonging to $\mathcal A_\alpha(\mu)$, then $\mathfrak m^{\mu,c}_{1,\beta}$ maps $L^\alpha(w^{\frac{\alpha}{s}}d\mu)$ into weak $L^s(wd\mu)$, i.e., there exists a constant $C>0$ such that for any element $f$ in $L^\alpha(wd\mu)$ and any real $a>0$,

$$\left(\int_{E_a}w(x)d\mu(x)\right)^{\frac{1}{s}}\leq \frac{C}{a}\left(\int_{\mathbb R^d}\left|f(x)\right|^\alpha w(x)^{\frac{\alpha}{s}}d\mu(x)\right)^{\frac{1}{\alpha}},$$
where $E_a=\left\{x\in\mathbb R^d/\mathfrak m^{\mu,c}_{1,\beta}f(x)>a\right\}$.
\end{thm}
In the case where $\beta=\infty$, Orobitg and P\'erez have shown (see Theorem 3.1 in \cite{OP}) that Theorem \ref{MuWee} remains true without the hypothesis that $\mu$ is the Lebesgue measure on $\mathbb R^d$.

It is clear that $\mathfrak m^{\mu,c}_{q,\beta}\leq \mathfrak m^{\mu}_{q,\beta}$. Notice that when $\mu$ is doubling, there exists a constant $C>0$ such that $\mathfrak m^{\mu}_{q,\beta}\leq C\mathfrak m^{\mu,c}_{q,\beta}$ and therefore Theorem \ref{MuWee} remains true with $\mathfrak m^{\mu}_{1,\beta}$ in place of $\mathfrak m^{\mu,c}_{1,\beta}$.

We recall that $\mu$ is said to be doubling if there exists a constant $C>0$ such that for any cube $Q$
$$\mu(2Q)\leq C\mu(Q),$$
where $2Q$ is the cube with same center as $Q$ but with side length twice that of $Q$.

However in the case where $\mu$ is nondoubling the picture is quite different. 

Actually we have the following result which is an extension to fractional maximal function of the theorem of Sj\"ogren in \cite{Sj}.
\begin{thm}\label{thm1.3}
\begin{enumerate}
\item[a)] Let $d=1$, $1\leq \alpha<\beta$, $\frac{1}{s}=\frac{1}{\alpha}-\frac{1}{\beta}$. If $w$ is a weight belonging to $\mathcal A_\alpha(\mu)$, then $\mathfrak m^{\mu}_{1,\beta}$ maps $L^\alpha(w^{\frac{\alpha}{s}}d\mu)$ into weak $L^s(wd\mu)$. 
\item[b)] For $d=2$ there is a measure $\mu$ for which $\mathfrak m^{\mu}_{1,\beta}$ does not map $L^1(\mu)$ into weak $L^{\beta'}(\mu)$
\end{enumerate} 
\end{thm}
The Gaussian measure $d\mu(x)=e^{-\frac{\left|x\right|^2}{2}}dx$ satifies the requirement of point $b)$. The proof is essentially the same as that of Theorem 3.2.6 of \cite{AABB}. To see it, just take $f=\chi_{B((0,a+1),\delta)}$ where $a>0$ large and $0<\delta<\frac{1}{a}$. There exist two positives absolutes constants $C_1$ and $C_2$ such that 
$$\left\|f\right\|_1\leq C_1\delta^2e^{-\frac{(a+1)^2}{2}}$$
 and 
$$\mu(B((s,a+1),1))^{\frac{1}{\beta}-1}\int_{B((s,a+1),1)}fd\mu>(\frac{e^{-\frac{a^2}{2}}}{a\sqrt{a}})^{\frac{1}{\beta}-1}C_2\delta^2e^{-\frac{(a+1)^2}{2}}$$
for all $\left|s\right|<1$. Thus, $\mathfrak m^{\mu}_{1,\beta}f$ is at least $\alpha:=(\frac{e^{-\frac{a^2}{2}}}{a\sqrt{a}})^{\frac{1}{\beta}-1}C_2\delta^2e^{-\frac{(a+1)^2}{2}}$ in the set $\left\{(x,y):\left|x\right|<1,a<y<a+2\right\}$, whose $\mu$-measure is at least $\frac{C_3}{a}e^{-\frac{a^2}{2}}$. Hence
$$\frac{\left\|f\right\|_1}{\alpha}\leq\frac{C_1}{C_2}(\frac{1}{\sqrt{a}})^{1-\frac{1}{\beta}}C^{\frac{1}{\beta}-1}_3\mu\left(\left\{(x,y)\in\mathbb R^2:\mathfrak m^{\mu}_{1,\beta}f(x,y)>\alpha\right\}\right)^{1-\frac{1}{\beta}},$$
which allows us to conclude. Notice that point a) can not be obtained by the classical method developed by Sj\"ogren if $w\neq 1$. In fact, the maximal operator is defined with the measure $\mu$, while the weak estimate is stated between two weighted spaces. Examples of such weights after Theorem \ref{thm1.3}, other than $w=1$, can be found in \cite{OP}.

Point $a)$ is a particular case of one of our main result (Theorem \ref{contopmax}). This result gives weighted norm inequalities for $\mathfrak m^{\mu}_{q,\beta}f$ when $f$ belongs to a family of Banach spaces $X^{q,p,\alpha}$ ($1\leq q\leq \alpha\leq p$) (see Section 2 for the definition) larger than that of Lebesgue spaces.

As a consequence, we obtain similar weighted norm inequalities for a class of convolution type  operators defined by
$$Kf(x)=\int_{\mathbb R}k(x-y)f(y)d\mu(y).$$
The class of convolution operators under consideration, contains the classical Riesz potential $I_\gamma$ ($0<\gamma<1$) defined by 
\begin{equation*}
I_{\gamma}f(x)=\int_{\mathbb R}\frac{f(x-y)dy}{\left|y\right|^{1-\gamma}},
\end{equation*}
whenever the integral exists.



The remaining of the paper is organized as follows:
In the $2^{nd}$ Section, we give prerequisites for the spaces $X^{q,p,\alpha}$. In Section 3, using an adapted Besicovitch covering lemma we prove the continuity of the fractional maximal operator $\mathfrak m^{\mu}_{q,\beta}$ 
between $X^{q,p,\alpha}$ spaces and weak Lebesgue spaces. In Section 4, we show that the weak-norm of the operator $K$
is controlled by the one of the fractional maximal operator $\mathfrak m^{\mu}_{q,\beta}$, a result which
allows us to deduce from the continuity of the fractional maximal operator, a weighted
inequality for some convolution type operators. In Section 5 we apply our results to an absolute
continuous measure with respect to the Lebesgue measure, and we show that some
classical results of harmonic analysis are special cases of ours.
Throughout the whole paper, $C$ denotes a positive constant which is independent
of the main parameters, but may change from line to line. The value of constants
with subscript as $C_0$ does not change in different occurrences. Whenever no precision
is given, $p,q$ and $\alpha$ will always stand for elements of $\left[1,+\infty\right]$ such that $q\leq\alpha\leq p$.
\section{Prerequisites on $X^{q,p,\alpha}$ spaces}
In all what follows, we denote by $\mu$ a  positive and non-atomic Radon measure on $\mathbb R$ satisfying 
\begin{equation}
\mu\left(\left(-\infty, a\right)\right)=\mu\left(\left[a, +\infty\right)\right)=\infty ,a\in\mathbb R.\label{mesure}
\end{equation}
We denote by $L^{0}:=L^{0}(\mathbb R,\mu)$ the complex vector space of equivalence
classes (modulo equality $\mu$- almost everywhere) of $\mu$ measurable
complex valued functions on $\mathbb R$.

We fix $x_0\in\mathbb R$. It is possible to associate to any real number $r>0$, a partition of $\mathbb R$ into intervals $I^{i}_{r}=\left[a^{r}_{i},a^{r}_{i+1}\right)$, $i\in\mathbb Z$ such that $a_0^r=x_0$ and $\mu\left(I_i^r\right)=r$ for all $i\in\mathbb{Z}$. For $1\leq p,q,\alpha\leq\infty$, put
\begin{equation*}
\|f\|_{q ,p,\alpha} = \sup_{r>0} r^{\frac{1}{\alpha} - \frac{1}{q}}\  _r\|f\|_{q,p}.\label{pqalpha}
\end{equation*}
where for $r>0$,
\begin{equation*}
_r \|f\|_{q, p} = \left\{
\begin{array}{lll}
\left[\sum \limits_{i \in \mathbb{Z}} \left(\|f\chi_{I_i^r}\|_q\right)^p%
\right]^{\frac{1}{p}} &\text{ if }&p<\infty  \\
\sup \limits_{i \in \mathbb{Z}} \|f\chi_{I_i^r}\|_q & \text{ if }
&p = \infty
\end{array}
\right.. \label{rpq}
\end{equation*}
$\left\|\cdot\right\|_{p}$ stands for the usual Lebesgue norm in $L^{p}:=L^{p}(\mathbb R,\mu)$ and $\chi_{I^{r}_{i}}$ is the characteristic function of the interval $I^{r}_{i}$. We consider in this paper, the spaces 
$$ X^{q, p}=\left\{f\in L^0/ _1 \|f\|_{q, p}<\infty\right\}$$
and
$$ \ X^{q,p,\alpha}=\{f\in L^0/\|f\|_{q,p,\alpha}<\infty \},$$
which are respectively the analogue of the Wiener amalgam space $(L^{q},\ell^{p})$ (see \cite{H,W}), and the space of integrable fractional mean functions $(L^{q},\ell^{p})^{\alpha}$ introduced by Fofana in \cite{Fo2}. Notice that the space $X^{q,p,\alpha}$ does not depend on the choice of $x_0\in\mathbb R$. The basic properties of these spaces are recapitulated in the following proposition.
\begin{prop}[\cite{AF}]Let $1\leq p,q,\alpha\leq\infty$. Then
\begin{enumerate}
\item   $\left(X^{q,p},\left\|\cdot\right\|_{q,p}\right)$ is a complex Banach space and
$\left(X^{q,p,\alpha},\left\|\cdot\right\|_{q,p,\alpha}\right)$ is a complete
normed subspace of $X^{q,p}$.
\item The space $X^{q,p,\alpha}$ is nontrivial if and only if $q\leq\alpha\leq p$.
\item If  $\alpha\in\left\{p,q\right\}$ then $X^{q,p,\alpha}=L^{\alpha}$.
\item If  $ q < \alpha < p$ then there exists $C>0$ depending only on $ p,q $ and $\alpha$ such that for all $f\in L^0$
\begin{equation*}
\|f\|_{q , p , \alpha} \leq C \|f\|_{\alpha ,\infty}^{\ast}, 
\end{equation*}
where
\begin{equation*}
\left\|f\right\|_{\alpha,\infty}^{\ast}=\sup_{\lambda>0}\lambda\mu\left(\left\{x\in\mathbb R:\left|f(x)\right|>\lambda\right\}\right)^{\frac{1}{\alpha}}.\label{lebesgue-faible}
\end{equation*}
\end{enumerate}
\end{prop}
Recall that the subspace of $L^{0}$ consisting of elements $f$ satisfying $\left\|f\right\|_{\alpha,\infty}^{\ast}<\infty$ is the weak Lebesgue space denoted by $L^{\alpha,\infty}$.

Weighted norm inequalities have been established between weak Lebesgue spaces and $(L^{q},\ell^{p})^{\alpha}$ spaces for fractional maximal operator  and  Riesz potential  in the context of Euclidean space \cite{Fo1}. 
 The following theorem has been proved in \cite{Fo1} in the Euclidean space, and generalized in \cite{FFK} in the context of spaces of homogeneous type.
\begin{thm}\label{theoB}Let $1\leq q\leq\alpha\leq p$ with  $0<\frac{1}{s}=\frac{1}{\alpha}-\frac{1}{\beta}$, $ q\leq q_{1}\leq\alpha_{1}\leq p_{1}$ with $0<\frac{1}{t}=\frac{1}{q_{1}}-\frac{1}{\beta}\leq\frac{1}{p_{1}}$,
and a weight $v$, satisfying 
\begin{equation*}
\sup_{\mathbb R\supset I \text{ interval } }\left|I\right|^{\frac{1}{\beta}-\frac{1}{q}}\left\|v\chi_{I}\right\|_{t}\left\|v^{-1}\chi_{I}\right\|_{1/(\frac{1}{q}-\frac{1}{q_{1}})}<\infty.
\end{equation*}
Then there exists a constant $C > 0$ such that
\begin{equation*}
\left(\int_{\left\{x\in\mathbb R:\mathfrak m^{dx}_{1,\beta}f(x)>\lambda\right\}}v(y)^{t}dy\right)^{\frac{1}{t}}\leq C\lambda^{-1}\left\|fv\right\|_{q_{1},p_{1},\alpha_{1}}\left(\lambda^{-1}\left\|f\right\|_{q,\infty,\alpha}\right)^{s\left(\frac{1}{q_{1}}-\frac{1}{\alpha_{1}}\right)}
\end{equation*}
for any real $\lambda>0$ and any Lebesgue measurable function $f$ on $\mathbb R$.
\end{thm}
This result, coupled with a weighted weak norm inequality between a fractional integral $I_{\gamma}$ and an appropriate fractional maximal operator  allows us to obtain a similar result for $I_{\gamma}$.


\section{Continuity of the maximal operator $\mathfrak m^{\mu}_{q,\beta}$} 
 Let $1\leq q\leq\beta\leq\infty$. For $f\in L^{q}_{loc}$, the fractional maximal function $\mathfrak m^{\mu}_{q,\beta}f$ is defined on the real line $\mathbb R$ as in (\ref{opmax}), with the cubes replaced by
 intervals in $\mathbb R$ containing $x$ and of finite measure. The first main result of this paper, which is the analogue of Theorem \ref{theoB}, can be stated as follows.

\begin{thm}\label{contopmax} Let $1 \leq q \leq \alpha \leq \beta$ and $q \leq q_1 \leq \alpha_1 \leq p_1$ such that $0< \frac{1}{q_1} - \frac{1}{\beta} =\frac{1}{\theta} \leq \frac{1}{p_1}$. 
Let $v$ be a weight on $\mathbb R$ for which there exists $C_{0}>0$ such that we have
  \begin{equation*}
\left\{\begin{array}{lll}
(\frac{1}{\mu(I)}\; \int_{I} \;v^{\theta} d\mu)^{\frac{1}{\theta}}(\frac{1}{\mu(I)}\int_{I}v^{-\frac{q_1 q}{q_1 -q}} d\mu)^{\frac{
q_1 - q}{q_1 q}}\leq C_{0}&\text{ if }&q < q_1\\
(\frac{1}{\mu(I)} \int_{I}v^{\theta} d\mu)^{\frac{1}{\theta}}(\|v^{-1} \chi_I\|_{\infty})\leq C_{0}&\text{ if }&q = q_1\end{array}\right.\label{cond1}
\end{equation*}
for any bounded interval $I\subset\mathbb R$. Then there exists $C_1 > 0$ such that for any $f\in L^0$
and $\lambda > 0$ we have :
\begin{enumerate}
\item $
\left(\int_{E_{\lambda}}v^{\theta} d\mu\right)^{\frac{1}{\theta}}\leq
C_{1}\lambda^{-1}\left(\int_{\mathbb R}\left|fv\right|^{q_{1}}d\mu\right)^{\frac{1}{q_{1}}}$,
\item if $\alpha<\beta$ then
\begin{equation}
\left(\int_{E_{\lambda}}v^{\theta} d\mu\right)^{\frac{1}{\theta}}\leq
C_{1}\left(\lambda^{-1}\|fv\|_{q_1 ,p_1 ,\alpha_1}\right)(\lambda^{-1}\;\|f\|_{q ,\infty
,\alpha})^{s(\frac{1}{q_1} - \frac{1}{\alpha_1})},\label{nouveau1}
\end{equation}
where $E_{\lambda} = \left\{x \in \mathbb{R} : \mathfrak m^{\mu}_{q,\beta}f(x)>\lambda \right\}$ and $\frac{1}{s} = \frac{1}{\alpha} - \frac{1}{\beta}$.
\end{enumerate}
\end{thm}
For the proof of this theorem, we will need the following covering lemma proved in \cite{AF}. A similar result to Theorem \ref{contopmax} is also given there in the context of measure.
\begin{lem}\label{recov}
Assume that :
\begin{itemize}
\item $\mathcal{F}$ is a family of left-closed intervals in $\mathbb R$, such that 
$$\sup\left\{\mu(F): F\in\mathcal{F}\right\}<\infty,$$
\item for any $F=\left[a,\: b\right)\in\mathcal{F}$, $c_F$ is a point of the open interval $\left(a, b\right)$ satisfying  $\mu\left(\left[a,c_F\right)\right)=\mu\left(\left[c_F, b\right)\right)=\frac{1}{2}\mu(F),$
\item $\mathfrak C=\{c_F : F\in\mathcal{F}\}.$
\end{itemize}
If $\left[A,B\right)$ is an interval of $\mathbb{R}$ such that $\mu\left(\left[A, B\right)\right)<\infty$,
then there exists a sequence $(F_i)_{i\in I}$ of elements of $\mathcal{F}$ satisfying  
$$\mathfrak C\cap\left[A, B\right)\subset\cup_{i\in I}F_i\text{ and }\sum_{i\in I}\chi_{F_i}\leq 5.$$
\end{lem}
\proof[Proof of Theorem \ref{contopmax}]
Let $f\in L^0$ and $\lambda>0$. For any $x\in E_{\lambda}$, there exists an interval $G_x\subset\mathbb{R}$ which contain $x$ and satisfies
 \begin{equation}
 \mu(G_x)^{\frac{1}{\beta}-\frac{1}{q}}\|f\chi_{G_x}\|_{q}> \lambda, \label{max1}
 \end{equation}
and two reals $a_x$ and $b_x$ satisfying 
\begin{equation}
a_x < x < b_x \text{ and }\mu(\left[a_x,x\right))=\mu(\left[x,b_x\right))=\mu(G_x).\label{div}
\end{equation}
We pose $I_x=\left[a_x,b_x\right)$.

Fix $x'\in E_{\lambda}$ and $R>\mu(G_{x'})>0$. There exist two reals $A$ and $B$ such that $\mu(\left[A,x'\right))=\mu(\left[x',B\right))=R$. Let $E_{\lambda,R}=\left\{x\in E_{\lambda}:G_{x}\subset\left[A,B\right)\right\}$. We have 
$$\sup\left\{\mu(I_{x}):x\in E_{\lambda,R}\right\}<\infty$$
since $\mu(\left[A,B\right))=2R$. Thus there exists a sub-family $(I_i)_{i \geq 1}$ of $ \{ I_x : x \in E_{\lambda}\}$ such that 
$$\sum_{i \geq 1}\chi_{I_i}\leq 5\text{ and }E_{\lambda , R}\subset\cup_{i \geq 1}I_i,$$
according to Lemma \ref{recov}. It follows that 
\begin{equation*}
\int_{E_{\lambda,R}}v^{\theta}d\mu \leq \int_{\cup_{i\geq1}I_i}v^{\theta} d\mu\leq \sum_{i \geq
1}\int_{I_i}v^{\theta} d\mu\leq \left[\sum_{i
\geq 1}(\int_{I_i}v^{\theta} d\mu)^{\frac{p_1}{\theta}}\right]^{\frac{\theta}{p_1}},
\end{equation*}
where we use for the last inequality, the fact that $\frac{p_1}{\theta}\leq 1$. Hence,
\begin{equation*}
\left(\int_{E_{\lambda,R}}v^{\theta} d\mu \right)^{\frac{p_1}{\theta}}
\leq \sum_{i \geq 1}\left(\int_{I_i}v^{\theta} d\mu\right)^{\frac{p_1}{\theta}}.
\end{equation*}
The choice of $G_x$ for $x\in E_{\lambda}$ and the construction of the $I_{x}$'s ensure that 
\begin{equation*}
\mu(I_x)^{\frac{1}{\beta}-\frac{1}{q}}\|f\chi_{I_x}\|_{q}= 2^{\frac{1}{\beta}-\frac{1}{q}}\mu(G_x)^{\frac{1}{\beta}-\frac{1}{q}}\|f\chi_{I_x}\|_{q}> 2^{\frac{1}{\beta}-\frac{1}{q}
}\mu(G_x)^{\frac{1}{\beta}-\frac{1}{q}}\|f\chi_{G_x}\|_{q}>2^{\frac{1}{\beta}-\frac{1}{q}}\lambda.
\end{equation*}
Thus $2^{\frac{1}{q}-\frac{1}{\beta}}\lambda^{-1}\mu(I_i)^{\frac{1}{\beta}-\frac{1}{q}}\|f\chi_{I_i}\|_{q}> 1$ for all $i\geq 1$, so that 
\begin{equation*}
\left(\int_{E_{\lambda,R}}v^{\theta} d\mu\right)^{\frac{p_1}{\theta}}\leq \underset{i \geq1}{\sum}\left(\int_{I_i}v^{\theta} d\mu\right)^{\frac{p_1}{\theta}}\left[2^{\frac{1}{q}-\frac{1}{\beta}}\lambda^{-1}\mu(I_i)^{\frac{1}{\beta}-\frac{1}{q}}\|f\chi_{I_i}\|_{q}
\right]^{p_1}.
\end{equation*}
The above relation becomes 
\begin{equation}
\left(\int_{E_{\lambda,R}}v^{\theta} d\mu\right)^{\frac{p_1}{\theta}}\leq \underset{i \geq1}{\sum}\left[\|v\chi_{I_i}\|_{\theta}\lambda^{-1}\|f\chi_{I_i}v\|_{q_{1}}\|v^{-1}
\chi_{I_i}\|_{\frac{q_{1} q}{q_{1} -q}}2^{\frac{1}{q} - \frac{1}{\beta}}\mu(I_i)^{\frac{1}{\beta} - \frac{1}{q}}\right]^{p_{1}},\label{qq}
\end{equation}
according to H\"older inequality. Consequently, given the assumption on $v$, there exists a real number $C$ not depending on $f$ and $\lambda$ such that 
\begin{equation}
\left(\int_{E_{\lambda,R}}v^{\theta}d\mu\right)^{\frac{p_1}{\theta}}\leq(C\lambda^{-1})^{p_{1}}\sum_{i \geq 1}\left[\|fv\chi_{I_i}\|_{q_1}\right]^{p_1}.\label{int}
\end{equation}
\begin{enumerate}
\item Since $q_{1}\leq p_{1}$, we have 
$$\sum_{i \geq 1}\left[\|fv\chi_{I_i}\|_{q_1}\right]^{p_1}\leq\left(\sum_{i \geq 1}\|fv\chi_{I_i}\|^{q_{1}}_{q_1}\right)^{\frac{p_{1}}{q_1}}\leq5^{\frac{p_{1}}{q_1}}\left\|fv\right\|^{p_{1}}_{q_{1}}.$$
 Taking the above relation in Estimate (\ref{int}) yields, 
\begin{equation}
\left(\int_{E_{\lambda,R}}v^{\theta}d\mu\right)^{\frac{1}{\theta}}\leq C \lambda^{-1}\left\|fv\right\|_{q_{1}}.\label{eqa1}
\end{equation}
\item We suppose that $\alpha<\beta$. If $\left\|f\right\|_{q,\infty,\alpha}=\infty$, then there is nothing to prove. We assume that $\left\|f\right\|_{q,\infty,\alpha}<\infty $. For all $x\in E_{\lambda}$, we have 
\begin{eqnarray*}
\mu(G_x)^{\frac{1}{q}-\frac{1}{\beta}} &\leq& \lambda^{-1}\|f\chi_{G_x}\|_{q}\leq\lambda^{-1}\sum_{k\in\mathbb Z}\left\|f\chi_{G_{x}\cap I^{\mu(G_{x})}_{k}}\right\|_{q}\leq2\lambda^{-1}\ _{\mu(G_{x})}\|f\|_{q,\infty}\\
&\leq&2\lambda^{-1}\mu(G_{x})^{\frac{1}{q}-\frac{1}{\alpha}}\|f\|_{q,\infty,\alpha}.
\end{eqnarray*}
Therefore 
$$\mu(I_{x})=2\mu(G_{x})\leq2\left(2\lambda^{-1}\|f\|_{q,\infty,\alpha}\right)^{s},$$
where $\frac{1}{s} = \frac{1}{\alpha} - \frac{1}{\beta}$. It follows that 
$r=\sup\left\{\mu(I_{x}):x\in E_{\lambda}\right\} \leq2\left(2\lambda^{-1}\|f\|_{q,\infty,\alpha}\right)^{s}.$ For any integer $j$, we put 
$$M_j=\left\{ i \geq 1 : \mu(I_i \cap I_{j}^r) >0\right\}\text{ and }M^{0}_{j}=\left\{i \geq 1 : \mu(I_i \cap I_{j}^r) = \mu(I_i)\right\}.$$
We have 
\begin{equation*}
\sum_{i\in M_{j}}\left(\int_{I_{i}\cap I^{r}_{j}}\left|fv\right|^{q_{1}}d\mu\right)^{\frac{p_{1}}{q_{1}}}\leq\left(\int_{\mathbb R}\left|fv\right|^{q_{1}}\sum_{M_{j}}\chi_{I_{i}\cap I^{r}_{j}}\right)^{\frac{p_{1}}{q_{1}}}\leq 5^{\frac{p_{1}}{q_{1}}}\left\|fv\chi_{I^{r}_{j}}\right\|^{p_{1}}_{q_{1}},
\end{equation*}
for all $j\in\mathbb Z$, so that 
\begin{equation}
\sum_{j\in\mathbb Z}\left\|fv\chi_{I^{r}_{j}}\right\|^{p_{1}}_{q_{1}}\geq 5^{-\frac{p_{1}}{q_{1}}}\sum_{j\in\mathbb Z}\sum_{i\in M_{j}}\left(\int_{I_{i}\cap I^{r}_{j}}\left|fv\right|^{q_{1}}d\mu\right)^{\frac{p_{1}}{q_{1}}}.\label{pp}
\end{equation}
The right hand side of (\ref{pp}) is greater than or equal to 
$$5^{-\frac{p_{1}}{q_{1}}}\left(\sum_{i\in \cup_{j\in\mathbb Z}M^{0}_{j}}\left(\int_{I_{i}}\left|fv\right|^{q_{1}}d\mu\right)^{\frac{p_{1}}{q_{1}}}+\sum_{j\in\mathbb Z}\sum_{i\in M_{j}\setminus M^{0}_{j}}\left(\int_{I_{i}\cap I^{r}_{j}}\left|fv\right|^{q_{1}}d\mu\right)^{\frac{p_{1}}{q_{1}}}\right),$$
since $\int_{I_{i}\cap I^{r}_{j}}\left|fv\right|^{q_{1}}d\mu=\int_{I_{i}}\left|fv\right|^{q_{1}}d\mu$ for all $i\in M^{0}_{j}$. Hence 
\begin{equation*}
\sum_{j\in\mathbb Z}\left\|fv\chi_{I^{r}_{j}}\right\|^{p_{1}}_{q_{1}}\geq C\sum_{i\geq 1}\left\|fv\chi_{I_{i}}\right\|^{p_{1}}_{q_{1}},
\end{equation*}
 since for $i\notin \cup_{j\in\mathbb Z}M^{0}_{j}$ (which implies that $i\in\cup_{j}(M_{j}\setminus M^{0}_{j})$), there exists a unique integer $j_{i}$ such that $I_{i}=(I_{i}\cap I^{r}_{j_{i}})\cup(I_{i}\cap I^{r}_{j_{i}+1})$. It follows therefore that 
\begin{equation}
\sum_{i \geq 1}(\|fv\chi_{I_i}\|_{q_1})^{p_1} \leq C\sum_{j \in \mathbb{Z}}\left(\left\|fv\chi_{I_{j}^r}\right\|_{q_1}\right)^{p_1}\leq C\left\|fv\right\|^{p_{1}}_{q_{1},p_{1},\alpha_{1}}r^{\frac{p_{1}}{\alpha_{1}}-\frac{p_{1}}{q_{1}}}.\label{0,5}
\end{equation}
Taking (\ref{int}) in (\ref{0,5}) with the control on $r$, we obtain
\begin{equation*}
\left(\int_{E_{\lambda,R}}v^{\theta} d\mu \right)^{\frac{1}{\theta}}\leq C
\lambda^{-1}\left\|fv\right\|_{q_{1},p_{1},\alpha_{1}}(\lambda^{-1}\left\|f\right\|_{q,\infty,\alpha})^{s (\frac{1}{q_1}
- \frac{1}{\alpha_1})}.
 \end{equation*}
The result follows by letting $R$ goes to infinity, since $v$ is positive.
\end{enumerate}
\epf

\begin{cor}\label{cor1}
Let $1 \leq q \leq \alpha < \beta$ with $0<\frac{1}{q} - \frac{1}{\beta}  \leq \frac{1}{p}\leq\frac{1}{\alpha}$. Then there exists $C>0$ such that 
\begin{equation*}
\left\|\mathfrak m^{\mu}_{q,\beta}f\right\|^{\ast}_{s}
\leq C\;\|f\|_{q,p,\alpha},\ f\in L^0
\end{equation*}
with $\frac{1}{s} = \frac{1}{\alpha}-\frac{1}{\beta} > 0.$
\end{cor}
\proof
We just have to take in Theorem \ref{contopmax} $v \equiv 1$, $\alpha_{1}=\alpha$ and $q_1 = q $
\epf

The next result is obtained by interpolation and the continuously embedding of $L^{\alpha}$ into $X^{q , p , \alpha}$.
\begin{cor}
Let $1 \leq q <\alpha <\beta$ with $\frac{1}{s} = \frac{1}{\alpha} - \frac{1}{\beta}$. Then there exists $C>0$ which depends only on $q,\beta$ and $\alpha$ such that 
\begin{equation*}
\left\|\mathfrak m^{\mu}_{q ,\beta}f\right\|_{s}\leq C\left\|f\right\|_{\alpha},\ f\in L^{\alpha}.
\end{equation*}
\end{cor}
\section{A control of some class of convolution type operators }
In this section, we fix $1\leq q\leq\beta\leq\infty$ and a positive even function $k$ on $\mathbb R$, which is non increasing on $\mathbb R_+$ and satisfies 
 $$\sup_{x\in\mathbb R}\left\|k(x-\cdot)\right\|^{\ast}_{\eta,\infty}<\infty\text{ and }\sup_{y\in\mathbb R}\left\|k(\cdot-y)\right\|^{\ast}_{\eta,\infty}<\infty,$$
 where  $k(\cdot-y):x\mapsto k(x-y)$, $k(x-\cdot):y\mapsto k(x-y)$ and $\frac{1}{\eta}=1-\frac{1}{\beta}$. We define the operator $K$ by 
\begin{equation}
Kf(x)=\int_{\mathbb R}k(x-y)f(y)d\mu(y),\label{opp}
\end{equation}
whenever the expression on the right hand side has a sense. Our first result in this section is the proof of a weak-norm inequality between $Kf$ and a suitable maximal function $\mathfrak m^{\mu}_{q,\beta}f$ for $f\in X^{q,p,\alpha}$. This result generalizes the analogue in the Euclidean space (see Theorem 1 in \cite{MW}).

\begin{thm}\label{contmaxpo}
We suppose that: 
\begin{enumerate}
\item there exists $C_0$ such that 
\begin{equation}
\overline{\displaystyle \lim_{r\rightarrow \infty}}\left(\sup_{t\in\mathbb R} \frac{\mu([t,
t+r])}{\mu([0, r])}\right) \leq C_0,\label{cnd1}
\end{equation}
and
\begin{equation}
\overline{\displaystyle \lim_{r\rightarrow \infty}}\left(\sup_{t\in\mathbb R} \frac{\mu([t-r,
t])}{\mu([-r, 0])} \right)\leq C_0,
\end{equation}
\item $k$ is lower semi-continuous function.
\item $\rho$ is a positive Borel measure on $\mathbb R$ satisfying an $\mathcal A_\infty$ condition, i.e., for all $\varepsilon > 0$ there exists $\delta > 0$ such that for any interval $I\subset\mathbb R$ we have 
$$\mu(E) \leq \delta \mu(I) \; \Rightarrow
\rho(E) \leq \varepsilon \;\rho(I),\ E\subset I.$$\label{gg}
\end{enumerate}
Then there exists $C>0$ such that
\begin{equation*}
\sup\limits_{\lambda >  0} \lambda^{\kappa}\; \rho\left(\{x \in \mathbb{R} : |Kf(x)| > \lambda \}\right) \leq C \sup\limits_{\lambda >  0}
\lambda^{\kappa}\; \rho\left(\{x \in \mathbb{R}:\mathfrak m^{\mu}_{q,\beta}f(x) > \lambda
\}\right),
\end{equation*}
for all $\kappa>0\text{ and }f\in X^{q,p,\alpha}$, where  $\frac{1}{p}=\frac{1}{q}-\frac{1}{\beta}$.
\end{thm}
For the proof of this result, we need the following lemmas.
\begin{lem}\label{lemme1} Let $p,q$ and $\beta$ be as in Theorem \ref{contmaxpo}. There exist two constants $B > 0$ and $D_1>0$ such that if: 
\begin{enumerate}[(i)]
\item $a, b$ and $c$ are real numbers satisfying $a > 0, b \geq B, c
> 0,$
\item $f$ is a positive element of $X^{q, p,\alpha}$,
\item $I = \left(x_1,  x_2\right)$ is an interval of $\mathbb{R}$ satisfying 
$\mu(I) < \infty$ and $Kf(x_j) \leq a$ for $j = 1,2,$
\end{enumerate}
then
$$\mu(\{x \in I: Kf(x) > ab  \text{ and }\mathfrak m^{\mu}_{q,\beta}f(x) \leq ac \}) \leq
D_1(\frac{c}{b})^p \mu(I).$$
\end{lem}
\proof
We consider a positive element $f$ of $X^{q,p,\alpha}$, real numbers $a > 0 , b > 0$ and $c > 0$, and an interval  $I\subset\mathbb{R}$ as in the statement. We assume that the set $\{x \in I:\mathfrak m^{\mu}_{q,\beta}f(x) \leq ac\}\neq\emptyset$, since otherwise there is nothing to prove.
Let $g = f\chi_I$ and $h = f-g.$
\begin{enumerate}
\item We have $g\in L^q$, since $\mu(I) < \infty$ and $f\in X^{q, p,\alpha}$. Besides that, since 
  $\frac{1}{\eta} + \frac{1}{q} - 1 = \frac{1}{p}$, it comes from Lemma 15.3  of \cite{FS} that $Kg \in L^{p,\infty}$
 and
 \begin{equation}
\left\|Kg\right\|_{p,\infty}^{\ast}\leq C\left\|k\right\|_{\eta,\infty}^{\ast} \left\|g\right\|_q.\label{ajout1}
\end{equation}
Putting together the definition of $L^{p,\infty}$ and Estimate (\ref{ajout1}), we obtain 
\begin{equation*}
\mu(\{x \in \mathbb{R}: |Kg| > \frac{ab}{2} \}) \leq \left[\frac{2C}{ab}\left\|k\right\|_{\eta,\infty}^{\ast} \left\|g\right\|_q \right]^p=D_1(\frac{1}{ab}\left\|g\right\|_q)^p,
\end{equation*}
with $D_1 = (2C\|k\|_{\eta,\infty}^{\ast})^p$.

However, for all $t\in \{x \in I: \mathfrak m^{\mu}_{q,\beta}f(x) \leq ac\}$, we have 
$$\left\|g\right\|_q = \left\|f\chi_I\right\|_q \leq \mu(I)^{\frac{1}{p}}\mathfrak m^{\mu}_{q,\beta}f(t) \leq ac \mu(I)^{\frac{1}{p}}.$$
 Therefore 
\begin{equation*}
\mu(\{x \in \mathbb{R}: |Kg| > \frac{ab}{2} \}) \leq D_1(\frac{c}{b})^p \mu(I).
\end{equation*}
\item Let $x\in I$. We have 
$$Kh(x) = \int_{-\infty}^{x_1} k(x-y) f(y)d\mu(y) \; + \;
\int_{x_2}^{+\infty} k(x-y) f(y) d\mu(y).$$
For the first term on the right hand side, we have 
$$\int_{-\infty}^{x_1} k(x-y) f(y)d\mu(y) \leq \int_{-\infty}^{x_1} k(x_1-y)
f(y)d\mu(y) \leq Kf(x_1) \leq a,$$
thanks to the non-increasing property of $k$ on $\mathbb R_+$.
Similarly,
$$\int_{x_2}^{+\infty} k(x-y) f(y) d\mu(y) \leq a.$$
Hence,
$$Kh(x) \leq 2a\ \ ,x\in I.$$
\item Let $B > 4$. If $b \geq B$ then we have $Kh(x) < \frac{ab}{2}$ for all $x\in I$, and consequently,
$$\begin{aligned}
&\mu(\{x \in I : Kf(x) > ab \text{ and } \mathfrak m^{\mu}_{q,\beta}f(x) \leq ac \})\\
&\ \ \ \ \ \ \ \leq \mu(\{x \in I : Kg(x) > \frac{ab}{2} \}) + \mu(\{x \in I : Kh(x) > \frac{ab}{2} \})\leq D_1(\frac{c}{b})^p \; \mu(I)
\end{aligned}$$
\end{enumerate}
Which achieves the proof of the lemma.
\epf

\begin{lem}\label{lemme2}Let $q,\beta$ and $\eta$ be as in Theorem \ref{contmaxpo}. There exists a constant $D_2>0$ such that if: 
\begin{enumerate}
\item $x_1, x_2, y_1 \text{ and } y_2$ are real numbers satisfying $y_1 < x_1 < x_2 < y_2 $  and $\mu([y_1 , x_1]) = \mu([x_2  , y_2]) \geq \mu([x_1  , x_2]),$
\item $f$ is a positive element of $L^0$ supported in the closed interval $\left[x_1 , x_2\right],$
\end{enumerate}
then
\begin{equation*}
\left\{\begin{array}{lll} Kf(x)\leq D_2 \left(\frac{\mu([x_2,x])}{\mu([0,x-x_2])}\right)^{\frac{1}{\eta}} \mathfrak m^{\mu}_{q,\beta}f(x)&\text{  if }&x >y_{2} ,\\
Kf(x) \leq D_2 \left(\dfrac{\mu([x , x_1])}{\mu([x-x_1 , 0])}\right)^{\frac{1}{\eta}}\mathfrak m^{\mu}_{q,\beta}f(x)&\text{ if }&x < y_1.
\end{array}
\right.
\end{equation*}
\end{lem}
\proof
Let $x_1,x_2,y_1,y_2 \in\mathbb R$, and $f\in L^0$ be as in the statements of the lemma.

Fix $x\in \left(y_2 , \infty\right]$. We have
$$Kf(x) = \int^{x_{2}}_{x_{1}} k(x-y)f(y) d\mu(y)\leq\int^{x_{2}}_{x_{1}} k(x-x_{2})f(y) d\mu(y)\leq k(x-x_2) \|f\|_q \mu(\left[x_{1},x_{2}\right])^{\frac{
1}{q^{\prime}}},$$
thanks to H\"older inequality. It follows that 
\begin{equation}
Kf(x) \leq k(x-x_2) \mu(\left[x_{1},x_{2}\right])^{\frac{1}{q'}}\left(\mu([x_1,x])\right)^{\frac{1}{q}-\frac{1}{\beta}} \mathfrak m^{\mu}_{q,\beta}f(x), \label{equa1}
\end{equation}
since $\mathrm{supp}f\subset\left[x_{1},x_{2}\right]\subset\left[x_{1},x\right]$. Put $s=\frac{\mu([x_2 , y_2 ])}{\mu(\left[x_{1},x_{2}\right])}$. We have 
$$\mu(\left[x_{1},x_{2}\right]) = \frac{1}{s}\mu([x_2 , y_2 ]) \leq \frac{1}{s}\mu([x_2 , x ]),$$
 which implies that 
$$\mu([x_1 , x]) = \mu(\left[x_{1},x_{2}\right]) + \mu([x_2 , x]) \leq (1 + \frac{1}{s}) \mu([x_2
, x]).$$
 So,
\begin{equation}
\mu(\left[x_{1},x_{2}\right])^{\frac{1}{q'}} (\mu([x_1,x]))^{\frac{1}{q}-\frac{1}{\beta}}  \leq
  2^{\frac{1}{q}-\frac{1}{\beta}} \mu([x_2 , x])^{\frac{1}{\eta}}.\label{extimat}
\end{equation}
Since $k\in L^{\eta,\infty}$ and it is non-increasing on $\mathbb R_{+}$, we have  
$$\mu([0,r]) \leq \mu(\{x \in \mathbb{R} : k(x) > \frac{k(r)}{2}\}) \leq
\left(\frac{2\|k\|_{\eta, \infty}^{\ast}}{k(r)}\right)^{\eta}, r > 0,$$
so that $k(r) \leq \frac{2\|k\|_{\eta, \infty}^{\ast}}{(\mu([0, r]))^{\frac{1}{\eta}}}$. If we take $r=(x-x_{2})$ and consider Relation (\ref{extimat}), then Inequality (\ref{equa1}) becomes 
$$Kf(x) \leq  D_2 \left[\frac{\mu([x_2, x])}{\mu([0, x-x_2])}\right]^{\frac{1}{\eta}
}\mathfrak m^{\mu}_{q,\beta}f(x).$$
An analogue argument allows us to show that for $x\in(-\infty,y_1)$, we have 
$Kf(x) \leq D_2 \left[\frac{\mu([x, x_1])}{\mu([x-x_1, 0])}\right]^{\frac{1}{\eta}}\mathfrak m^{\mu}_{q,\beta}f(x)$
\epf

\proof[Proof of Theorem \ref{contmaxpo}]
\begin{enumerate}
\item Let $f$ be a positive element of $X^{q,p,\alpha }$.

$1^{rst} $ case: we suppose that $f$ has compact support embedded in an interval $J=[x_1, x_2]$ of $\mathbb{R}$. For $\lambda > 0$, put 
$$F_{\lambda} = \{x\in \mathbb{R} : |Kf(x)| > \lambda \}.$$
We have $\mu(F_{\lambda} )<\infty$, since $Kf\in L^{p,\infty}$. We also have that $Kf$ is lower semi-continuous according to Proposition 2.3.2 of \cite{AH}, since by hypothesis $k$ is lower semi-continuous. It follows that $F_{\lambda}$ is open in $\mathbb R$, and therefore can be wrote as disjoint union of open intervals; i.e., $F_{\lambda} = \cup_{m\in M} I_m$, where the $I_{m}$'s are disjoint open intervals of $\mathbb R$ and $M$ is a countable set. Let $m\in M$ and $c>0$ a real number. We have 
$$\mu\left(\{x \in I_m : Kf(x) > \lambda B\text{ and }\mathfrak m^{\mu}_{q,\beta}f(x) \leq \lambda c\}\right) \leq D_1 (\frac{c}{B})^p \mu(I_m)$$
according to Lemma \ref{lemme1}, where the constants $B$ and $D_{1}$ are those appearing in the statement of the lemma. It follows from the assumptions on the measure $\rho$ that there exists
$\delta > 0$ such that  for all $m\in M$ 
$$\rho\left(\{x \in I_{m} : Kf(x) >\lambda B\text{ and } \mathfrak m^{\mu}_{q,\beta}f(x) \leq \lambda c\}\right) \leq \frac{B^{-\kappa}}{2} \rho\left(I_{m}\right),$$
whenever $0 < c\leq \delta$. By the fact that $I_m$'s are disjoint, we deduce that 
 $$\rho\left(\{x \in \mathbb R : Kf(x) >\lambda B\text{ and } \mathfrak m^{\mu}_{q,\beta}f(x) \leq \lambda c\}\right) \leq \frac{B^{-\kappa}}{2} \rho\left(F_{\lambda}\right).$$
But it is easy to see that
$$\left\{x\in\mathbb R:Kf(x)>\lambda B\right\}\subset\left\{x\in\mathbb R:\mathfrak m^{\mu}_{q,\beta}f>\lambda c\right\}\cup\left\{x\in\mathbb R:Kf(x) >\lambda B\text{ and }\mathfrak m^{\mu}_{q,\beta}f(x)\leq\lambda c\right\} .$$
Therefore
$$\begin{aligned}
&\rho\left(\{x \in \mathbb{R} : Kf(x)> \lambda B \}\right) \\
&\ \ \ \ \ \ \ \ \ \ \ \leq
\rho\left(\{x \in \mathbb{R} :\mathfrak m^{\mu}_{q,\beta}f(x) > \lambda c \} \right) +
\rho\left(\{x \in \mathbb{R} :Kf(x) >\lambda B\text{ and }\mathfrak m^{\mu}_{q,\beta}f(x)\leq \lambda c \} \right),
\end{aligned}$$
and consequently,
\begin{equation}
\rho\left(\{x \in \mathbb{R} : Kf(x)> \lambda B \}\right) \leq \rho\left(\{x \in \mathbb{R}:\mathfrak m^{\mu}_{q,\beta}f(x) >\lambda c \} \right) +
\frac{B^{-\kappa}}{2} \rho\left(F_{\lambda}\right).\label{etoile}
\end{equation}
We put $L=\left[y_1 , y_2\right]\text{ with }y_{1}<x_{1}<x_{2}<y_{2}\text{ and }\mu(\left[y_{1},x_{1}\right])=\mu(\left[x_{2},y_{2}\right])=\mu(\left[x_{1},x_{2}\right])$.
For all $x\in\mathbb{R}\setminus L$, we have
$$\left\{\begin{array}{lll}Kf(x) \leq D_2 \left(\frac{\mu([x_2 , x])}{\mu([0 , x-x_2])}\right)^{\frac{1}{\eta}}\mathfrak m^{\mu}_{q,\beta}f(f)(x)&\text{ if }&x > y_2\\
Kf(x) \leq D_2 \left(\frac{\mu([x , x_1])}{\mu([x-x_1 , 0])}\right)^{\frac{1}{\eta}} \mathfrak m^{\mu}_{q,\beta}f(x)&\text{ if }&x < y_1.
\end{array}\right.,$$
according to Lemma \ref{lemme2}. Thus, by choosing the interval $[x_1,x_2]$ large enough, we have 
$$ Kf(x) \leq D_2 (C_0)^{\frac{1}{\eta}}\mathfrak m^{\mu}_{q,\beta}f(x),$$
for all $x\in\mathbb R\setminus L$, thanks to Relation (\ref{cnd1}).
We choose $c = \inf \{\delta , \frac{1}{D_2 (C_0)^{\frac{1}{\eta}}} \}$. For all $x \in \mathbb{R}\setminus L$ such that $Kf(x) > \lambda $, we have 
$$ D_2 (C_0)^{\frac{1}{\eta}}\mathfrak m^{\mu}_{q,\beta}f(x) > \lambda,$$
that is,
\begin{equation*}
\mathfrak m^{\mu}_{q,\beta}f(x) >\frac{1}{D_2 (C_0)^{\frac{1}{\eta}}} \lambda\geq\lambda c.
\end{equation*}
Thus 
$$\{x \in \mathbb{R} \setminus L : Kf(x) > \lambda \} \subset \{x \in
\mathbb{R} : \mathfrak m^{\mu}_{q,\beta}f(x) > \lambda c \}.$$
Therefore, we have 
 \begin{eqnarray*}
\rho\left(\{x \in \mathbb{R}: Kf(x) > \lambda B \}\right) &\leq& \rho\left(\{x \in \mathbb{R} :\mathfrak m^{\mu}_{q,\beta}f(x) >  \lambda c \} \right)\\
&+& \frac{B^{-\kappa}}{2} \left[\rho\left(\{x \in \mathbb{R} \setminus L :
Kf(x) > \; \lambda \} \right)+\rho\left(\{x \in L : Kf(x) >  \lambda \}\right)\right] \\
&\leq &2 \rho\left(\{x \in \mathbb{R} :\mathfrak m^{\mu}_{q,\beta}f(x) > \lambda c \}
\right) + \frac{B^{-\kappa}}{2} \rho\left(\{x \in L : Kf(x) > \lambda \}
\right),
\end{eqnarray*}
according to Relation (\ref{etoile}). Multiplying both sides of the above inequality
by $\lambda^{\kappa}$ and taking the $\sup$ for all $0<\lambda<N$, we obtain 
\begin{eqnarray*}
\sup _{0 <\lambda< N} \lambda^{\kappa} \;\rho\left(\{x \in \mathbb{R} : Kf(x) > \lambda B \}\right)  &\leq &2 \sup _{0<
\lambda < N}\lambda^{\kappa} \;\rho\left(\{x \in \mathbb{R} :\mathfrak m^{\mu}_{q,\beta}f(x) >\lambda c \} \right) \\
&+& \frac{1}{2}\sup_{0<\lambda < N} B^{-\kappa}\lambda^{\kappa} \;\rho\left(\{x \in L : Kf(x) >\lambda \} \right).
\end{eqnarray*}
The above relation can also be written as
\begin{eqnarray*}
B^{-\kappa}\sup_{0<\lambda <NB}  \lambda^{\kappa}\;
\rho\left(\{x \in \mathbb{R} : Kf(x) > \lambda \}\right)  &\leq& 2 \sup_{0 < \lambda < Nc} c^{-\kappa} \lambda^{\kappa}\;
\rho\left(\{x \in \mathbb{R} :\mathfrak m^{\mu}_{q,\beta}f(x) > \lambda \} \right)\\
&+& \frac{B^{-\kappa}}{2}\sup_{0 < \lambda < N}
\lambda^{\kappa} \;\rho\left(\{x \in L : Kf(x) >\lambda \} \right),
\end{eqnarray*}
that is 
\begin{equation*}
\frac{B^{-\kappa}}{2}\sup_{0<\lambda <NB}  \lambda^{\kappa}\;
\rho\left(\{x \in \mathbb{R} : Kf(x) > \lambda \}\right)  \leq 2 \sup_{0 < \lambda < Nc} c^{-\kappa} \lambda^{\kappa}\;
\rho\left(\{x \in \mathbb{R} :\mathfrak m^{\mu}_{q,\beta}f(x) > \lambda \} \right).
\end{equation*}
Letting $N$ goes to infinity, we have
$$\sup \limits_{\lambda  >  0} \lambda^{\kappa}\; \rho\left(\{x \in
\mathbb{R} : Kf(x) > \lambda \}\right) \leq C \sup \limits_{\lambda > 0} \lambda^{\kappa} \;\rho\left(\{x \in \mathbb{R} :\mathfrak m^{\mu}_{q,\beta}f(x) >
\lambda \} \right)$$

$2^{nd}$ case : we suppose that the support of $f$ is not necessarily compact.

For all integers $n \geq 1$, we put $f_n = f\chi_{[-n,n]}$. We have that for all positive integers $n$, $f_n\in X^{q, p,\alpha}$ and has compact support. It follows from the first case that 
\begin{equation}
\sup _{\lambda  > 0} \lambda^{\kappa} \;\rho\left(\{x \in
\mathbb{R} : Kf_n(x) > \lambda \}\right) \leq C \sup_{\lambda
>0} \lambda^{\kappa} \;\rho\left(\{x \in \mathbb{R} :\mathfrak m^{\mu}_{q,\beta}f_{n}(x) > \lambda \} \right).
\end{equation}
Since the sequence $(Kf_n)_{n \geq 1}$ is an increasing sequence which converges to $Kf$ and the sequence $(\mathfrak m^{\mu}_{q,\beta}f_{n})_{n \geq 1}$ is bounded above by $\mathfrak m^{\mu}_{q,\beta}f$, the monotone convergence theorem leads to the expected result.\label{partie1}
\item  For an arbitrary element $f$ of $X^{q,p,\alpha}$, let $f_1 = \sup(f, 0)$ and $f_2 = \sup(-f,
0)$. Both $f_1$ and $f_2$ are positive elements of $X^{q, p, \alpha}$ with $f = f_1 - f_2$. We also have 
$$\{x \in \mathbb{R} : |Kf(x)| >  \lambda \} \subset \{x \in \mathbb{R}
: Kf_1(x) >  \frac{\lambda}{2} \} \cup \{x \in \mathbb{R} : Kf_2(x)> \frac{\lambda}{2} \},$$
so that the result follows from part \ref{partie1} and the facts that $\mathfrak m^{\mu}_{q,\beta}f_{1} \leq \mathfrak m^{\mu}_{q,\beta}f$ and $\mathfrak m^{\mu}_{q,\beta}f_{2} \leq \mathfrak m^{\mu}_{q,\beta}f$.
\end{enumerate}
\epf

Before the next result, we recall the definition of condition $\mathcal A_p(\mu)$. Let $1 < p <\infty$. A weight $w$ satisfies the condition (or belongs to the class) $\mathcal A_p(\mu)$ if
$$\sup_{\mathbb R\supset I:\text{ interval}}\left(\frac{1}{\mu(I)}\int_{I}w(x)^{-\frac{1}{p-1}}d\mu(x)\right)^{p-1}\left(\frac{1}{\mu(I)}\int_{I}w(x)d\mu(x)\right)< +\infty.$$
We put $\mathcal A_\infty(\mu)=\cup_{p>1}A_p(\mu)$.
\begin{prop}\label{contoppo}We suppose that
\begin{itemize}
\item there exists $C_0>0$ such that 
\begin{equation} \overline{ \lim_{r\rightarrow
\infty}}\left(\sup_{t\in\mathbb R}\frac{\mu([t, t + r])}{\mu([0, r])}\right) \leq C_0 \text{ and }\overline{ \lim_{r\rightarrow \infty}}\left(\sup_{t\in\mathbb R}\frac{\mu([t - r,
t])}{\mu([-r, 0])}\right) \leq C_0,\label{controlemesure}
\end{equation}
\item $k$ is lower semi-continuous,
\item $1\leq q,\alpha,\beta,q_{1},\alpha_{1},p_{1}\leq\infty$ satisfy
$$\left\{\begin{array}{lll}0\leq\frac{1}{\beta}\leq\frac{1}{\alpha}\leq\frac{1}{q}\leq1&&\\
q\leq q_{1}\leq\alpha_{1}\leq p_{1}&\text{ with }&0<\frac{1}{\theta}=\frac{1}{q_{1}}-\frac{1}{\beta}
\end{array}\right. $$
\item $v$ is a positive measurable function on $\mathbb R$ and there exists a real constant
$M > 0$ such that for all intervals $I$ of $\mathbb R$ we have 
\begin{equation}
\left\{\begin{array}{lll}
\left(\frac{1}{\mu(I)} \int_{I} v^{\theta} d\mu\right)^{\frac{1}{\theta}}\left(\frac{1}{\mu(I)} \int_{I} v^{-\frac{q_1 q}{q_1 - q}} d\mu \right)^{%
\frac{q_1 -q}{q_1 q}} \leq M &\text{  if  }&q < q_1,\\
\left(\frac{1}{\mu(I)} \int_{I} v^{\theta} d\mu\right)^{\frac{1}{\theta}}
\|v^{-1} \chi_I\|_{\infty} \leq M&\text{  if  }&q = q_1.
\end{array}\right.
\end{equation}
\end{itemize}
Then there exists $C>0$ such that for any positive element $f$ of $L^0$ and $0<\lambda<\infty $
\begin{enumerate}
\item $\left(\int_{F_{\lambda}}v^{\theta}d\mu\right)^{\frac{1}{\theta}}\leq C\left(\lambda^{-1}\int_{\mathbb R}\left|fv\right|^{q_{1}}d\mu\right)^{\frac{1}{q_{1}}}$, where  $F_{\lambda} = \{x \in \mathbb{R} : Kf(x)> \lambda\}$. 
\item If $\frac{1}{\beta}<\frac{1}{\alpha}$, then we have 
\begin{equation}
\left(\int \limits_{F_{\lambda}} v^{\theta} d\mu \right)^{\frac{1}{\theta}}
\leq c \lambda^{-1} \|fv\|_{q_1 , p_1 ,\alpha_1} \; \left(\lambda^{-1}
\|f\|_{q, \infty,\alpha}\right)^{s(\frac{1}{q_1} - \frac{1}{\alpha_1})}
\end{equation}
where $\frac{1}{s} = \frac{1}{\alpha} - \frac{1}{\beta}$.
\end{enumerate}
\end{prop}
\proof
Put $w = v^{\theta}$. We have $w\in\mathcal A_r(\mu)$, with $r= 1+\theta(\frac{1}{q} -
\frac{1}{q_1})$. Therefore, the measure $\rho$ such that $d\rho(x) = w(x)d\mu(x)$ satisfies Condition \ref{gg}) of Theorem \ref{contmaxpo}, according to Lemma 2.3 of \cite{OP}. It follows that
\begin{equation}
\sup _{\lambda  > 0} \lambda^{\kappa} \;\rho\left(\{x \in
\mathbb{R} : |Kf(x)| >  \lambda \}\right) \leq C \sup _{\lambda >  0} \lambda^{\kappa} \;\rho\left(\{x \in \mathbb{R} :\mathfrak m^{\mu}_{q,\beta}f(x) >\lambda \} \right),\ \ \kappa>0 .
\end{equation}
Since $F_{\lambda} = \{x \in \mathbb{R} : Kf(x)  > \lambda\}$ and
$E_{\lambda} = \{x \in \mathbb{R} :\mathfrak m^{\mu}_{q,\beta}f(x) > \lambda \}$, it comes that 
\begin{equation}
\sup_{\lambda > 0} \lambda^\kappa  \int_{F_{\lambda}} v^{\theta}
d\mu(x) \leq C_1 \sup_{\lambda > 0} \lambda^\kappa\int_{E_{\lambda}} v^{\theta} d\mu(x),\ \ \kappa>0.\label{ineq3}
 \end{equation}
The required result follows from (\ref{ineq3}) and Theorem \ref{contopmax}.
\epf

The next result is an immediate consequence of Proposition \ref{contoppo}.
\begin{cor}\label{corollaire3.5}We suppose that
$$\max\left\{\overline{ \lim_{r\rightarrow
\infty}}\left(\sup_{t\in\mathbb R}\frac{\mu([t, t + r])}{\mu([0, r])}\right),\overline{ \lim_{r\rightarrow \infty}}\left(\sup_{t\in\mathbb R}\frac{\mu([t - r,
t])}{\mu([-r, 0])}\right)\right\}<\infty.$$
\begin{enumerate}
\item If $1\leq q,\alpha,\beta,q_{1},\alpha_{1},p_{1}\leq\infty$ satisfy
\begin{itemize}
\item $1\leq q\leq\alpha\leq\beta$ with $0<\frac{1}{s}=\frac{1}{\alpha}-\frac{1}{\beta}$,
\item $q\leq q_{1}\leq\alpha_{1}\leq p_{1}$ with $0<\frac{1}{\theta}=\frac{1}{q_{1}}-\frac{1}{\beta}\leq\frac{1}{p_{1}}$
\end{itemize}
then there exists $C>0$ such that for all positive elements $f$ of $L^{0}$ and all $\lambda>0$ we have 
\begin{equation}
\mu(F_{\lambda})^{\frac{1}{\theta}}\leq C \lambda^{-[1+s(\frac{1}{q_{1}}-\frac{1}{\alpha_{1}})]}\left\|f\right\|_{q_{1},p_{1},\alpha_{1}}\left\|f\right\|^{s(\frac{1}{q_{1}}-\frac{1}{\alpha_{1}})}_{q,\infty,\alpha},
\end{equation}
where $F_{\lambda} = \{x \in \mathbb{R} : Kf(x)> \lambda\}$ .\label{premier}
\item In particular if $1\leq q\leq\alpha<\beta$, $\frac{1}{s}=\frac{1}{\alpha}-\frac{1}{\beta}$ and $\frac{1}{\theta}=\frac{1}{q}-\frac{1}{\beta}\leq\frac{1}{p}\leq\frac{1}{\alpha}$ then there exists $C>0$ such that for all positive elements $f$ of $L^{0}$ we have
\begin{equation*}
\left\|Kf\right\|^{\ast}_{s,\infty}\leq C\left(\left\|f\right\|^{\frac{1}{s}}_{q,p,\alpha}\left\|f\right\|^{\frac{1}{\theta}-\frac{1}{s}}_{q,\infty,\alpha}\right)^{\theta}\leq C\left\|f\right\|_{q,p,\alpha}.
\end{equation*}
\end{enumerate}
\end{cor}
\proof
\begin{enumerate}
\item It is clear that $v\equiv 1$ satisfy the hypotheses of Proposition \ref{contoppo}. Consequently,
the assertion follows from it.
\item Consider positive elements $f\in L^{0}$ and $\lambda>0$. Taking $q_{1}=q$, $\alpha_{1}=\alpha$ and $p_{1}=p$ in  (\ref{premier}), we obtain 
\begin{equation*}
\mu(F_{\lambda})^{\frac{1}{\theta}}\leq C\lambda^{-s/\theta}\left\|f\right\|_{q,p,\alpha}\left\|f\right\|^{s(\frac{1}{q}-\frac{1}{\alpha})}_{q,\infty,\alpha},
\end{equation*}
that is 
\begin{equation*}
\lambda\mu(F_{\lambda})^{\frac{1}{s}}\leq C\left(\left\|f\right\|^{\frac{1}{s}}_{q,p,\alpha}\left\|f\right\|^{\frac{1}{q}-\frac{1}{\alpha}}_{q,\infty,\alpha}\right)^{\theta}.
\end{equation*}
It follows that for every positive element $f$ of $L^0$, we have 
\begin{equation*}
\left\|Kf\right\|^{\ast}_{s,\infty}\leq C\left(\left\|f\right\|^{\frac{1}{s}}_{q,p,\alpha}\left\|f\right\|^{\frac{1}{\theta}-\frac{1}{s}}_{q,\infty,\alpha}\right)^{\theta}\leq C\left\|f\right\|_{q,p,\alpha}.
\end{equation*}
The last inequality comes from the fact that $\left\|f\right\|_{q,\infty,\alpha}\leq\left\|f\right\|_{q,p,\alpha}$.
\end{enumerate}
\epf

A consequence of the above result is the following corollary.
\begin{cor}\label{corollaire3.6}Suppose that:
$$\max\left(\overline{ \lim_{r\rightarrow
\infty}}\left(\sup_{t\in\mathbb R}\frac{\mu([t, t + r])}{\mu([0, r])}\right),\overline{ \lim_{r\rightarrow \infty}}\left(\sup_{t\in\mathbb R}\frac{\mu([t - r,
t])}{\mu([-r, 0])}\right) \right)<\infty.$$
If $1<\alpha<\beta<\infty$ then there exists $C>0$ such that for all positive elements $f$ of $L^{0}$ we have
\begin{equation*}
\left\|Kf\right\|^{\ast}_{s,\infty}\leq C\left\|f\right\|^{\ast}_{\alpha,\infty},
\end{equation*}
where $\frac{1}{s}=\frac{1}{\alpha}-\frac{1}{\beta}$.
\end{cor}
\proof
Since $0<\frac{1}{\beta}<\frac{1}{\alpha}<1$ and $\frac{1}{\alpha}-\frac{1}{\beta}<\frac{1}{\alpha}$, we can fine $p,q>1$ satisfying $\frac{1}{\alpha}-\frac{1}{\beta}<\frac{1}{q}-\frac{1}{\beta}\leq\frac{1}{p}<\frac{1}{\alpha}$. Note that $1\leq q<\alpha<p$. Thus, thanks to Proposition \ref{contoppo}, there exists $C_{1}>0$ such that for all positive elements $f$ of $L^{0}$ we have 
\begin{equation*}
\left\|f\right\|_{q,\infty,\alpha}\leq\left\|f\right\|_{q,p,\alpha}\leq C_{1}\left\|f\right\|^{\ast}_{\alpha,\infty}.
\end{equation*}
Moreover, by Corollary \ref{corollaire3.5}, there exists $C_{2}>0$ such that for all positive elements $f$ of $L^{0}$ we have 
\begin{equation*}
\left\|Kf\right\|^{\ast}_{s,\infty}\leq C_{2}\left\|f\right\|^{\frac{\theta}{s}}_{q,p,\alpha}\left\|f\right\|^{1-\frac{\theta}{s}}_{q,\infty,\alpha}.
\end{equation*}
The combination of these two inequalities completes the proof.
\epf

\section{Application}

We assume that $0<a<\gamma<1\leq\alpha<\frac{1-a}{\gamma-a}$ and let $\frac{1}{s}=\frac{1}{\alpha}-\frac{\gamma-a}{1-a}$.
 \begin{enumerate}
\item Denote by $\mu$ the measure defined on $\mathbb R$ by $d\mu(x)=\left|x\right|^{-a}dx$. We have
%
\begin{equation*}
\sup_{t\in\mathbb R}\frac{\mu(\left[t,t+r\right])}{\mu(\left[0,r\right])}\leq 2,\ r>0
\end{equation*}
and
\begin{equation*}
\sup_{t\in\mathbb R}\frac{\mu(\left[t-r,t\right])}{\mu(\left[-r,0\right])}\leq 2,\ r>0.
\end{equation*}
\item Define $k$ on $\mathbb R$ by $k(x)=\left|x\right|^{\gamma-1}$. We remark that $k$ is positive, lower semi-continuous, even, decreasing on $\mathbb R_+$ and satisfies
\begin{equation*}
\sup_{y\in\mathbb R}\left\|k(\cdot-y)\right\|^{\ast}_{\eta,\infty}\leq 2^{1-\gamma}(\frac{2}{1-a})^{\frac{1}{\eta}},
\end{equation*}
where $\frac{1}{\eta}=\frac{1-\gamma}{1-a}=1-\frac{\gamma-a}{1-a}$.
\item Consider the Riesz potential $I_{\gamma}$ defined by 
$$I_{\gamma}f(x)=\int_{\mathbb R}\left|x-y\right|^{\gamma-1}f(y)dy.$$
It is clear that for every positive element $f$ of $L^{0}$, we have
$$I_{\gamma}f(x)=\int_{\mathbb R}k(x-y)F(y)d\mu(y)\equiv KF(x),\ x\in\mathbb R$$
where $F(y) = f(y)|y|^a$ . Therefore, applying Corollaries \ref{corollaire3.5} and \ref{corollaire3.6}, we obtain the following result.
 \begin{prop}\label{4.1}
 \begin{enumerate}
\item If $1\leq q\leq\alpha$ and $\frac{1}{\theta}=\frac{1}{q}-\frac{\gamma-a}{1-a}\leq\frac{1}{p}\leq\frac{1}{\alpha}$ then there exists $C>0$ such that 
\begin{equation*}
\left\|I_{\gamma}f\right\|^{\ast}_{s,\infty}\leq C\left\|F\right\|^{\frac{\theta}{s}}_{q,p,\alpha}\left\|F\right\|^{1-\frac{\theta}{s}}_{q,\infty,\alpha}.
\end{equation*}
with $\frac{1}{s}= \frac{1}{\alpha}-\frac{\gamma-a}{1-a} , F(x)=\left|x\right|^{a}f(x)$, $x\in\mathbb R$ and $f$ being a positive element of $ L^{0}$.
\item If $1<\alpha$, then there exists $C>0$ such that 
\begin{equation}
\left\|I_{\gamma}f\right\|^{\ast}_{s,\infty}\leq C\left\|F\right\|^{\ast}_{\alpha,\infty},\label{4.1.b}
\end{equation}
for all positive elements $f$ of $L^{0}$.
 \end{enumerate}
\end{prop}
\end{enumerate}
Suppose that $1 < \alpha$. Taking into consideration Marcinkiewicz interpolation theorem, Relation (\ref{4.1.b}) of Proposition \ref{4.1} allows us to obtain the following classical result (see \cite{SW}). There exists $C > 0$ such that
$$\left\|I_{\gamma}f\right\|_{s}\leq C\left\|F\right\|_{\alpha} ,\ f\in L^{\alpha}(\mathbb R,\mu);$$
that is 
$$\left[\int_{\mathbb R}\left(\left|x\right|^{-\frac{a}{s}}I_{\gamma}f(x)\right)^{s}dx\right]^{\frac{1}{s}}\leq C\left[\int_{\mathbb R}\left|f(x)\right|^{\alpha}\left|x\right|^{a\frac{\alpha}{\alpha'}}dx\right]^{\frac{1}{\alpha}},\ f\in L^{\alpha}(\mathbb R,\left|x\right|^{a\frac{\alpha}{\alpha'}}dx).$$


\begin{thebibliography}{MTW1}
\bibitem{AABB}Claire Anantharaman, Jean-Philippe Anker, Martine Babillot, Aline Bonami, Bruno Demange,
et al..  Th\'eor\`emes ergodiques des actions de groupes.  L'Enseignement Math\'ematique.  41, 2010,
Monographies de L'Enseignement Math\'ematiques, 978-2-940264-08-7."hal-00464094"
\bibitem{AF} A. Adama et I. Fofana, {Continuit\'e de l'op\'erateur maximal fractionnaire associ\'e \`a une mesure de Radon positive}, Annales Africaines de Math\'ematiques serie.1, vol.2 (2011) 129-137.
\bibitem{AH} David R. Adams and Lars Inge Hedberg, Function spaces and potential theory, Grundlehren der mathematischen Wissenschaften 314  A Series of Comprehensive Studies in Mathematics, Springer.
\bibitem{Fo2} I. Fofana, {\it \'Etude d'une classe d'espaces de fonctions contenant les espaces de Lorentz}, Afr. Mat., ({\bf 1}) (1988) 29-50.
\bibitem{Fo1} I. Fofana, {\it Continuit\'e de l'int\'egrale fractionnaire et espace $(L^{q},\ell^{p})^{\alpha}$}, C. R. Acad. Sci. Paris 308, S\'erie I,(1989) 525-527.
\bibitem{FFK}J. Feuto, I. Fofana and K. Koua, {\it Weighted norms inequalities for a maximal operator in some subspace of amalgams}. Canad. Math. Bull. 53 ({\bf 2}) 2010, 263-277.
\bibitem{FS}G. Folland and E. M. Stein, Estimates for the $\bar{\partial}_{b}$ Complex and Analysis on the Heisenberg Group, Comm. Pure Appl. Math, XXVII (1974), 429-522.
\bibitem{H} F. Holland, {\it Harmonic Analysis on amalgams of $L^{p}$ and $\ell^{q}$}, J. London Math. Soc. (2) {\bf 10} (1975), 295-305.
\bibitem{MW} B. Muckenhoupt and R. Wheeden, {\it Weighted norm inequalities for fractional integrals }. Trans. Amer.Math. Soc. 192(1974), 261-274.
\bibitem{OP} J. Orobitg, C. Pérez, {\it $A_p$ Weights for nondoubling measures in $\mathbb{R}^n$ and applications },
Trans. Am. Math. Soc,354(2002), 2013-2033.
\bibitem{PW}C. P\'rez and R.L. Wheeden, {\it Uncertainty principle estimates for vector fields}, J. Functional
Analysis 181 (2001), 146–188
\bibitem{Sj}P. Sj\"ogren,{\it A remark on the maximal function for measures in $\mathbb R^n$}, Amer. J. Math. {\bf 105}
(1983), 1231-1233.
\bibitem{SS}P. Sj\"ogren and F. Soria, {\it Sharp estimates for the non-centered maximal operator associated to Gaussian and other radial measures}, Advances in Mathematics 181 (2004) 251-275.
\bibitem{SW} E. M. Stein and G. Weiss, {\it Fractional integral on $n$-dimensional Euclidean space}, J. Math. rech. {\bf 7} (1958), 503-514.
\bibitem{W} N. Wiener, {\it On the representation of functions by trigonometrical integrals}, Math. Z.24 (1926), 575-616.
\end{thebibliography}
\end{document}